\documentclass[12pt,a4paper]{amsart}
\usepackage{amssymb}
\usepackage{enumerate}
\usepackage{natbib}

\newtheorem{thm}{Theorem}[section]
\newtheorem{prop}[thm]{Proposition}
\newtheorem{lemma}[thm]{Lemma}

\newtheorem{dfn}[thm]{Definition}

\newcommand{\pf}{\noindent{\em Proof. }}

\DeclareMathOperator{\Aut}{Aut}
\DeclareMathOperator{\id}{id}
\DeclareMathOperator{\mlt}{Mlt}

\newcommand{\F}{\mathbb F}

\begin{document}

\title{Direct construction of code loops}
\author[G.P. Nagy]{G\'abor P. Nagy}
\email{nagyg@math.u-szeged.hu}
\address{Bolyai Institute, University of Szeged,
Aradi v\'ertan\'uk tere 1, H-6720 Sze\-ged, Hungary}
\urladdr{http://www.math.u-szeged.hu/\~{}nagyg/}
\thanks{The author was supported by the FKFP grant $0063/2001$ of the
Hungarian Ministry for Education and the OTKA grants nos. F042959 and
T043758.}


\begin{abstract}
Code loops were introduced by R. L. Griess. \citet{Griess} and \citet{Hsu}
gave methods to construct the corresponding code loop from any given doubly
even binary code; both these methods used some kind of induction. In this
paper, we present a global construction of the loop, where we apply the
correspondance between the concepts of Moufang loops and groups with
triality.
\end{abstract}

\maketitle

\section{Introduction}

For a doubly even binary code, one can introduce three operation in a
combinatorial way. These operations are related by
polarization. \citet{Hsu} called vector spaces over $\F_2$ with such
operations \emph{symplectic cubic spaces.} In \cite{Griess}, R.~L.~Griess
has introduced the notion of \emph{code loops;} these are Moufang loops $L$
such that $L/A$ is an elementary Abelian $2$-group for some central
subgroup $A$ of order $2$. In code loops, the power, the commutator and the
associator maps define a symplectic cubic space on $L/A$. The fact that any
symplectic cubic space corresponds to a code loop was shown by
\citet{Griess}; he constructed the \emph{factor set} of the loop extension
by induction on the dimension. \citet{Hsu} gave more explicit method to
construct the loops as \emph{centrally twisted products.} However, also
Hsu's method uses some inductive argument, since he puts together the loop
from smaller parts using a clever product rule. 

In this paper, we present a \emph{global construction} for the
loop. Therefore, we use that Moufang loops can be equivalently given by the
specific group theoretical concept of groups with triality. Hence, starting
from the symplectic cubic space, we first define a group $G$ and
automorphisms £$\sigma, \rho$ of $G$. Then, we prove that $G$ is a group
with triality with respect to these automorphisms. Finally, we show that
the Moufang loop $L$ corresponding to $G$ is a code loop, which gives rise
to the symplectic cubic space we started with. This global approach made
possible to implement our method using the \citet{GAP4} computer algebra
package by \citet{LOOPS}. 

For completeness, we mention \citet{Kitazume} where the author gives
another explicit construction of the code loop for some special cases of
doubly even codes. Also, the work of \citet{CheinGoodaire} is important,
where they showed the reverse impliciation, namely they constructed doubly
even binary codes from a given symplectic cubic space. Unfortunately, this
code is by far not unique; the one in the general construction has
parameters $[s,n]$, where $n$ is the dimension of the symplectic cubic
space and $s=O(\frac{2}{3} n^3)$. Finally, we mention the paper
\citet{Petr}, where the author generalizes the notions of symplectic cubic
spaces and doubly even codes and the construction of Chein and Goodaire.

\section{Preliminaries}

The subspace $\mathcal C \leq \F_2^n$ is a \emph{doubly even binary code,}
if $4|w(x)$ for all $x\in \mathcal C$, where the \emph{weight} $w(x)$ of $x$
is the number of nonzero coordinates of $x$. We can identify $\mathcal C$
with a set of subsets of $\{1,\ldots,n\}$. In this manner,
$w(x)=|x|$. Moreover, if $\mathcal C$ is doubly even, we have $2|w(x\cap
y)$ and the functions $\sigma:\mathcal C \to \F_2$, $\kappa: \mathcal C
\times \mathcal C \to \F_2$, $\alpha:\mathcal C \times \mathcal C \times
\mathcal C \to \F_2$, 
\begin{eqnarray*}
\sigma(x)&=&\frac{1}{4} \, w(x) \pmod{2}, \\
\kappa(x,y) &=&\frac{1}{2} \, w(x\cap y) \pmod{2}, \\
\alpha(x,y,z) &=&w(x\cap y \cap z) \pmod{2}.
\end{eqnarray*}
are well defined. Clearly, $\kappa$ and $\alpha$ are alternating and they
satify the relations
\begin{eqnarray}
\sigma(x+y)&=&\sigma(x)+\sigma(y)+\kappa(x,y), \label{eq:sigmapol} \\
\kappa(x+y,z)&=&\kappa(x,z)+\kappa(y,z)+\alpha(x,y,z), \label{eq:kappapol} \\
\alpha(x+y,z,t)&=&\alpha(x,z,t)+\alpha(y,z,t). \label{eq:alphalin} 
\end{eqnarray}
With other words, $\alpha$ is trilinear and $\kappa$ and $\alpha$ is
obtained from $\sigma$ and $\kappa$ by polarization, respectively. 

\medskip

The set $L$ together with a binary operation $(x,y) \mapsto x\cdot y=xy$ is
a \emph{quasigroup} if the equation $xy=z$ can be uniquely solved if two of
the three indeterminants is given. A quasigroup with a neutral element is a
\emph{loop.} Every element $x$ of a loop $L$ determines \emph{left} and
\emph{right multiplication maps} $L_{x},R_{x}: L\to L$, $yL_x=xy$ and
$yR_{x}=yx$. (As the reader can see, we write group actions on the right
hand side.) The maps $R_x,L_x$ are clearly permutations of $L$. The
permutation group $\mlt(L)$ generated by all left and right multiplications
is called the \emph{multiplication group of} $L$. The permutations $L_{x,y}
= L_xL_yL_{yx}^{-1}$, $R_{x,y} = R_xR_yR_{xy}^{-1}$ are called \emph{inner
maps} of $L$.

\emph{Moufang loops} are defined by the \emph{Moufang identity}
\[x(y(xz))=((xy)x)z.\] 
Recall that Moufang loops are \emph{diassociative}, which means that every
subloop generated by two elements is a subgroup. In particular, the power
$x^n$, $n\in \mathbb N$, the inverse $x^{-1}$ and the commutator
$[x,y]=x^{-1}y^{-1}xy$ are well-defined for every $x,y\in L$.

The \emph{associator} of $x$, $y$, $z\in L$ is the element $(x,y,z)=(xy\cdot
z)^{-1}(x\cdot yz)$. The \emph{center} $Z(L)$ of $L$ consists of the
elements $z\in L$ satisfying $[x,z]=(x,y,z)=(x,z,y)=(z,x,y)=1$ for all $x,y
\in L$. One sees easily that every subloop $A\leq Z(L)$ is a \emph{normal
subloop} of $L$, that is, the factor loop $L/A$ is well defined. (See
\citet{Bruck} for basic concepts on quasigroups and Moufang loops.)

The Moufang loop $L$ is called a \emph{small Frattini Moufang loop,} if
$L/A$ is an elementary Abelian $p$-group for some $A\leq Z(L)$ with
$|A|=p$. By \citet{Hsu}, for $p>3$, any small Frattini Moufang loop is
associative. Assume $L$ to be a small Frattini Moufang $2$-loop which is
not an elementary Abelian $2$-group. Let us identify the vector space $V$
over $\F_2$ with the factor loop $L/A$ and the subloop $A$ with the field
$\F_2$ of order $2$. We can introduce the following operations on $V$:
\[ \sigma(xA) = x^2, \hskip 5mm
\kappa(xA,yA) = [x,y], \hskip 5mm
\alpha(xA,yA,zA) = (x,y,z). \]
Then, the equations \eqref{eq:sigmapol}, \eqref{eq:kappapol} and
\eqref{eq:alphalin} hold for $\sigma, \alpha, \kappa$. 

\medskip 

\begin{dfn}
Let $V$ be a vector space over $\mathbb \F_2$. Let $\sigma:V \to \F_2$,
$\kappa:V\times V \to \F_2$ and $\sigma:V \times V \times V \to \F_2$ be
maps satisfying \eqref{eq:sigmapol}, \eqref{eq:kappapol} and
\eqref{eq:alphalin}. Then, $(V, \sigma, \kappa, \alpha)$ is called a
\emph{symplectic cubic space.} 
\end{dfn}

We mention that \eqref{eq:sigmapol}, \eqref{eq:kappapol} and
\eqref{eq:alphalin} imply
\begin{eqnarray*}
\kappa(x,x)&=& 0, \\
\kappa(x,y)&=& \kappa(y,x), \\
\alpha(x,y,z) &=& \sigma(x+y+z) + \sigma(x+y) + \sigma(y+z) + \sigma(x+z) + 
\\
&& \sigma(x) + \sigma(y) + \sigma(z). 
\end{eqnarray*}
Hence, $\alpha$ is a trilinear alternating form on $V$. 

\medskip

Let $S_n$ be the symmetric group on $\{1, \ldots, n\}$. We have the
following definition due to \citet{Doro}.

\begin{dfn}
The pair $(G,S)$ is called a {\em group with triality}, if $G$ is
a group, $S \leq \Aut{G}$, $S=\langle \sigma,\rho \mid
\sigma^2=\rho^3=(\sigma\rho)^2=1 \rangle \cong S_3$, and for all
$g\in G$ the {\em triality identity}
\[ [g,\sigma]\,[g, \sigma]^\rho\, [g,\sigma]^{\rho^2}=1\]
holds.
\end{dfn}

The following equivalent formulation of the concept of a group with
trialitiy is well known. 

\begin{lemma}[Parker] \label{lm:parker}
Let $G$ be a group and let $\sigma_1,\sigma_2,\sigma_3$ be involutorial
automorphisms of $G$. Let us denote by $\mathcal C_i$ the conjugacy class
$\sigma_i^G \subseteq \Aut(G)$. Then, $(G,\langle \sigma_1, \sigma_2
\rangle)$ is a group with triality if and only if $(\tau_i\tau_j)^3=\id$
for all $\tau_i \in \mathcal C_i$, $\tau_j \in \mathcal C_j$, $i, j \in
\{1,2,3\}$, $i\neq j$.
\end{lemma}

The next lemma characterizes a special class of groups with triality.

\begin{lemma} \label{lm:trgreq}
Let $G$ be a group and let $S\leq \Aut(G)$ be isomorphic to the symmetric
group $S_3$ on $3$ elements. Let us denote by $\sigma_1,\sigma_2, \sigma_3$
the involutions of $S$ and put $H_i=C_G(\sigma_i)$ and $\mathcal
C_i=\sigma_i^G$. Then, if $H_i$ acts transitively on $\mathcal C_j$ for
some $i\neq j$, then $(G,S)$ is a group with triality. In particular, if
$|G:H_j| = |H_i:H_i \cap H_j| < \infty$ for $i\neq j$, then $(G,S)$ is a
group with triality.
\end{lemma}
\pf Choose arbitrary $\tau_i=\sigma_i^g \in \mathcal C_i$ and $\tau_j =
\sigma_j^h \in \mathcal C_j$ with $i\neq j$. By the assumption, there is $f
\in H_i$ such that $\sigma_j^f = \sigma_j^{hg^{-1}}$. Then, 
\[\tau_i\tau_j = (\sigma_i \sigma_j^{hg^{-1}})^g = (\sigma_i\sigma_j)^{fg},\]
which implies $(\tau_i\tau_j)^3=((\sigma_i\sigma_j)^3)^{fg}=\id$. By Lemma
\ref{lm:parker}, $(G,S)$ is a group with triality. 

\qed

\section{Constructing the group with triality}
\label{sec:grconstr}

Let $V=(V,\sigma,\kappa,\alpha)$ be a symplectic cubic space. Let us choose
a basis $\mathcal{B} = \{b_1, \ldots, b_n\}$ of $V$ and denote by
$\sigma_i$, $\kappa_{ij}$ and $\alpha_{ijk}$ the structure constants of $V$
with respect to $\mathcal B$. 

We define the group $G$ with gerenators $g_i, f_i, h_i$, $i\in \{1,\ldots,
n\}$, $u$ and $v$ by the following relations:
\begin{eqnarray}
&& g_i^2 = u^{\sigma_i}, \; f_i^2 = v^{\sigma_i}, \; h_i^2=u^2=v^2=1,
  \label{eq:G1} \\
&& [g_i,g_j] = u^{\kappa_{ij}}, \; [f_i,f_j] = v^{\kappa_{ij}},
  \label{eq:G2} \\ 
&& [g_i,f_j] = (uv)^{\kappa_{ij}} \, \prod_{k=1}^n h_k^{\alpha_{ijk}},
  \label{eq:G3} \\ 
&& [g_i,h_j] = u^{\delta_{ij}}, \; [f_i,h_j] = v^{\delta_{ij}},
  \label{eq:G4} \\ 
&& [h_i,h_j] = [g_i,u] = [f_i,u] = [h_i,u] = [g_i,v] = [f_i,v] = [h_i,v] = 
  1. \label{eq:G5} 
\end{eqnarray}

\begin{lemma}
The group $G$ is well defined. Any element of $G$ is of the form $g_1^{x_1}
\cdots g_n^{x_n} \, f_1^{y_1} \cdots f_n^{y_n} \, h_1^{z_1} \cdots h_n^{z_n}
\, u^{t_1}v^{t_2}$ with $x_i,y_i,z_i,t_i \in \mathbb Z_2$. In particular,
the order of $G$ is $2^{3n+2}$.
\end{lemma}
\pf In order to show that $G$ is well defined, we prove the following.
\begin{enumerate}[(i)]
\item $E=\langle f_1,\ldots,f_n, h_1, \ldots, h_n, v\rangle$ is an
extraspecial $2$-group of type $+$ and order $2^{2n+1}$.
\item The commutator and power relations for the $g_i$'s are consistent. In
  particular, modulo $\langle u \rangle$, the group $\langle g_1, \ldots,
  g_n \rangle$ is an elementary Abelian group of order $2^n$.
\item The map $\gamma_i$ induced on $E \times \langle u\rangle$ by $g_i$ is
  an automorphism.
\item The $\gamma_i$'s induce an automorphism group $A$ of $E \times
\mathbb Z_2$ which is an elementary Abelian $2$-group.
\end{enumerate}
On the one hand, in $E$, $\langle h_1,\ldots, h_n, v\rangle$ is a maximal
elementary Abelian $2$-group of order $2^{n+1}$. On the other hand, with
\[\tilde f_i = f_i \, h_i^{\sigma_i} \, \prod_{k>i} h_k^{\kappa_{ik}},\]
we have
\[ \begin{array}{rcl}
[ \tilde f_i, \tilde f_j ] &=& [f_i,f_j] 
   [h_i^{\sigma_i} \prod_{k>i}   h_k^{\kappa_{ik}}, f_j ] 
   [f_i, h_j^{\sigma_j} \prod_{k>j} h_k^{\kappa_{jk}} ] \\
&=& v^{\kappa{ij}} v^{\kappa{ij}} \\
&=&1 
   \end{array} \]
and
\[ \tilde f_i^2 = f_i^2 [f_i, h_i^{\sigma_i} \prod_{k>i} h_k^{\kappa_{ik}} ]
= v^{\sigma_i} v^{\sigma_i} =1.\]
This means that $\langle \tilde f_1, \ldots, \tilde f_n \rangle$ is
elementary Abelian. Finally, $[\tilde f_i, h_j]
=[f_i,h_j]=v^{\delta_{ij}}$, hence $E$ is as stated in (i).

Let $\kappa^*$ be the anternating bilinear form on $V$ with structure
constants $\kappa_{ij}$ and $q$ be the quadratic form obtained by the
quadratic extension of $\sigma$ with respect to $\kappa^*$. Then, $q$
determines a central extension of $V$, isomorphic to $\langle g_1, \ldots,
g_n, u \rangle$. This proves (ii).

(iii) follows from the fact that $\gamma_i$ preserves the relations of $E
\langle u\rangle$. Indeed, putting 
\[f'_j = \gamma_i(f_j) = f_j \, (\prod_k
h_k^{\alpha_{ijk}})(uv)^{\kappa_{ij}}, \hskip 1cm
h'_j=\gamma_i(h_j) = h_j u^{\delta_{ij}}, \]
we have
\begin{eqnarray*}
[f'_\ell, f'_j] &=& [f_\ell,f_j] [f_\ell, \prod_k h_k^{\alpha_{ijk}}] [f_j,
  \prod_k h_k^{\alpha_{i\ell k}}]{} \\
&=& v^{\kappa_{ij}} v^{\alpha_{ij\ell}} v^{\alpha_{i\ell j}} \\
&=& \gamma_i([f_\ell,f_j]), \\
{}[f'_\ell, h'_j] &=& [f_\ell, h'_j] [\prod_k h_k^{\alpha_{i\ell k}}, h'_j]{}
  \\
&=& [f_\ell, h_j] {} = \gamma_i([f_\ell, h_j]), \\
(f'_j)^2 &=& f_j^2 \, [f_j, \prod_k h_k^{\alpha_{ijk}}] = f_j^2 =
  \gamma_i(f_j^2). 
\end{eqnarray*}
To show (iv), we calculate the action of $\gamma_{i_1} \circ
\gamma_{i_2}$. One obtains
\begin{eqnarray*}
f_j &\mapsto& f_j \, (\prod_k h_k^{\alpha_{i_1jk}+\alpha_{i_2jk}}) \,
(uv)^{\kappa_{i_1j} + \kappa_{i_2j}} \, u^{\alpha_{i_1ji_2}}, \\
h_j &\mapsto& h_j \, u^{\delta_{i_1j} + \delta_{i_2j}}.
\end{eqnarray*}
This means $\gamma_{i_1} \circ \gamma_{i_2} = \gamma_{i_2} \circ
\gamma_{i_1}$ and $\gamma_i^2 = \id$, hence (iv) holds. The other
statements are trivial. \qed

\begin{lemma} \label{lm:welldef}
The maps
\begin{eqnarray}
\sigma&:& g_i \leftrightarrow f_i, \, h_i \mapsto h_i, \, u \leftrightarrow
v \label{eq:sigma} \\ 
\rho&:& g_i \mapsto f_i, \, f_i \mapsto (g_if_i)^{-1}, \, h_i \mapsto h_i,
\, u \mapsto v, \, v \mapsto uv \label{eq:rho} 
\end{eqnarray}
extend to automorphisms of $G$. Moreover, $\sigma^2=\rho^3=(\sigma\rho)^2
=\id$, hence, $S=\langle \sigma, \rho \rangle$ is isomorphic to the
symmetric group $S_3$ on $3$ elements.
\end{lemma}
\pf The fact that $\sigma$ is an involutorial automorphism is
trivial. Also, $\rho$ preserves the relations of $G$; we only present the
calculations in the two most complex cases. We first observe that $g_i$ and
$f_i$ commute, therefore $(g_if_i)^2 \in Z(G)$. Moreover,
$[g_i,f_j]=[f_i,g_j]$ commutes with $f_i,f_j,g_i,g_j$. Thus,
\begin{eqnarray*}
[\rho(f_i), \rho(f_j)] &=& [ (g_if_i)^{-1}, (g_jf_j)^{-1} ] \\
&=& [ g_if_i, g_jf_j ] \\
&=& [f_i,g_j] [f_i,f_j]^{g_jg_i} [g_i,g_j] [g_i,f_j]^{g_j} \\
&=& [f_i,f_j] [g_i,g_j] \\
&=& (uv)^{\kappa_{ij}} \\
&=& \rho( [f_i,f_j] ).
\end{eqnarray*}
Similarly, 
\begin{eqnarray*}
[\rho(g_i), \rho(f_j)] &=& [ f_i, (g_jf_j)^{-1} ] \\
&=& [ f_i, g_jf_j ] \\
&=& [f_i,f_j] [f_i,g_j]^{f_j} \\
&=& [f_i,f_j] [f_i,g_j] \\
&=& v^{\kappa_{ij}} \, (\prod_k h_k^{\alpha_{ijk}}) \, (uv)^{\kappa_{ij}}
\\
&=& \rho( [g_i,f_j] ).
\end{eqnarray*}
Finally, the relations for $\sigma$ and $\rho$ hold, since $\sigma^2$,
$\rho^3$ and $(\sigma\rho)^2$ leave the generators of $G$ invariant. \qed

\begin{lemma}
Let us define the subgroups
\begin{eqnarray*}
H_1 &=& \langle g_i, h_i, u \mid i=1,\ldots,n \rangle, \\
H_2 &=& \langle f_i, h_i, v \mid i=1,\ldots,n \rangle, \\
H_3 &=& \langle g_if_i, h_i, uv \mid i=1,\ldots,n \rangle 
\end{eqnarray*}
of $G$. Then, $H_3=C_G(\sigma)$, $H_1^\sigma=H_2$, $H_1^\rho=H_2$,
$H_2^\rho=H_3$. 
\end{lemma}
\pf Clearly, $H_3\leq C_G(\sigma)$. For the converse, let us write the
element $a\in G$ in the form 
\[ a= g_1^{x_1}f_1^{y_1} \cdots g_n^{x_n}f_n^{y_n} \, h_1^{z_1}
\cdots h_n^{z_n} \, u^{t_1}v^{t_2}. \]
One immediately has that $\sigma(a)=a$ only if $x_i=y_i$ and $t_1=t_2$,
that is, $a\in H_3$. This proves $H_3=C_G(\sigma)$, the rest is
trivial. \qed

\begin{prop}
Let $(V,\sigma,\kappa,\rho)$ be a symplectic cubic space and let us define
the group $G$ by \eqref{eq:G1}--\eqref{eq:G5}. Moreover, define the
automorphisms $\sigma$ and $\rho$ by \eqref{eq:sigma} and \eqref{eq:rho},
respectively. Then, $(G,\langle \sigma, \rho \rangle)$ is a group with
triality.
\end{prop}
\pf We have $H_1\cap H_2 = H_1 \cap H_3 = H_2 \cap H_3 = \langle h_q,
\ldots, h_n \rangle$, hence $|G:H_3| = |H_2:H_2 \cap H_3|=2^{n+1}$. By
Lemma \ref{lm:trgreq}, $(G,\langle \sigma, \rho \rangle)$ is a group with
triality. \qed

\section{Some properties of Moufang loops given by groups with triality}

In the $3$-net, the lines $X=a$, $Y=a$ and $XY=a^{-1}$ are permuted by
$\rho$. Indeed, let $\sigma_1, \sigma_2, \sigma_3$ denote the Bol
reflections with respect to the $Y$-axis $X=1$, $X$-axis $Y=1$ and
transversal line $XY=1$, respectively. Then,
\[ \begin{array}{l}
\sigma_1 : \left \{ \begin{array}{rcl}
X=a & \leftrightarrow & X=a^{-1}, \\
Y=a & \leftrightarrow & XY=a, \\ 
Y=a^{-1} & \leftrightarrow & XY=a^{-1}, 
\end{array} \right. \vspace{1mm} \\ 
\sigma_2 : \left \{ \begin{array}{rcl}
X=a & \leftrightarrow & XY=a, \\ 
Y=a & \leftrightarrow & Y=a^{-1}, \\ 
X=a^{-1} & \leftrightarrow & XY=a^{-1}, 
		    \end{array} \right. \vspace{1mm} \\
\sigma_3 : \left \{ \begin{array}{rcl}
X=a & \leftrightarrow & Y=a^{-1}, \\ 
Y=a & \leftrightarrow & X=a^{-1}, \\ 
XY=a & \leftrightarrow & XY=a^{-1}, 
		    \end{array} \right. 
\end{array}\]
and $\rho=\sigma_2\sigma_1$ acts as claimed. 

Moreover, in the coordinate loop, $ab=c$ holds if and only if the lines
$X=a$, $Y=b$, $XY=c$ are concurrent. Let us denote by $\tau_a$ the Bol
reflection with respect to the axis $X=a$. Then, $\tau_a^\rho$ and
$\tau_a^{\rho^2}$ are the Bol reflectiosn with respect to the lines $Y=a$
and $XY=a^{-1}$, respectively. The equation $ab=c$ holds if and only if
$\langle \tau_a, \tau_b^\rho, \tau_{c^{-1}}^{\rho^2} \rangle \cong S_3$,
that is, if and only if
\[ \tau_{c^{-1}} = (\tau_a \tau_b^\rho \tau_a)^\rho =
\tau_b^{\rho\tau_a\rho} . \]
Since $\tau_{c^{-1}} = \tau_1 \tau_c \tau_1$ and $\tau_1=\sigma_1=\sigma$,
we have
\begin{equation} \label{eq:netop}  
ab=c \; \Leftrightarrow \; \tau_c = \tau_b^{\rho \tau_a \rho \sigma}.
\end{equation}

Let $(G,S)$ be a group with triality, $S=\langle \sigma, \rho \mid
\sigma^2=\rho^3=(\sigma\rho)^2 =\id \rangle$. As before, we denote by
$\sigma=\sigma_1$, $\sigma_2$, $\sigma_3$ the involutions of $S$. The
conjugacy class $\sigma_i^G$ be $\mathcal C_i$. In \citet{HallNagy}, we
showed how to construct a (dual) $3$-net from $(G,S)$: $\mathcal C_1 \cup
\mathcal C_2 \cup \mathcal C_3$ are the lines and $\tau_i \in \mathcal
C_i$, $i=1,2,3$ are concurrent if and only if $\langle \tau_1, \tau_2,
\tau_3 \rangle \cong S_3$.

\begin{prop} \label{pr:mltgr}
Let $(G,S)$ be a group with triality and use the notation
$\sigma=\sigma_1$, $\rho$, $\mathcal C_1=\sigma^G$ as before. Let us define
the binary operation
\begin{equation} \label{eq:moufop} 
\alpha \circ \beta = \beta^{\rho\alpha\rho\sigma} =
\alpha^{\rho^{-1}\beta\rho^{-1}\sigma} 
\end{equation}
on $\mathcal C_1$. Then, 
\begin{enumerate}[(i)]
\item The operation is well defined, $\sigma$ is a two sided unit
  element. $(\mathcal C_1,\circ)$ is a isomorphic to the Moufang loop
  associated to $(G,S)$.
\item Put $\alpha=\sigma^g$ and $\gamma=[g,\sigma]^\rho$. Using the natural
  bijection between $\mathcal C_1 =\sigma^G$ and the set of right cosets
  $X=G/C_G(\sigma)$, the right action of $\gamma$ on $X$ is equivalent with
  the right multipliation $R_\alpha$ of the loop. Similarly, the right
  action of $\gamma^\rho$ on $X$ is equivalent with the left multiplication
  $L_\alpha$.
\item Let $N$ be the largest normal subgroup of $G$, contained in
  $C_G(\sigma)$. Then, the multiplication group of the Moufang loop
  associated to $(G,S)$ is a subgroup of $G/N$. 
\end{enumerate}
\end{prop}
\pf Clearly, \eqref{eq:netop} implies (i), and (ii) implies
(iii). Furthermore, on the one hand,
\[ \rho^{-1}\alpha\rho^{-1}\sigma = \rho^{-1} \sigma^g \sigma \rho =
   [g,\sigma]^\rho = \gamma.\]
On the other hand,
\[ \beta R_\alpha = \beta \circ \alpha = \beta^{\rho^{-1}\alpha\rho^{-1}
  \sigma} = \beta^\gamma.\]
hence (ii) follows. \qed

\medskip

The next lemma will make possible to calculate the structure constants of
small Frattini Moufang loops which are given by their groups with triality.

\begin{lemma} \label{lm:moufids}
Let $L$ be a Moufang loop satisfying $x^2 \in Z(L)$ for all $x\in L$. Then
\begin{enumerate}[(i)]
\item $[R_x,R_y]= R_{[x,y]}$, 
\item $[R_y,L_z]=R_{y^{-1},z} =L_{y,z^{-1}}$, 
\item $[[R_x,L_y],R_z] = R_{(x,y,z)}$ 
\end{enumerate}
hold for all $x,y,z \in L$. 
\end{lemma}
\pf On the one hand, by \citet[Theorem 1.1.6]{NSt_canad}, Moufang loops
satisfying $x^2 \in Z(L)$ for all $x\in L$ are conjugacy closed loops,
hence $R_y^{-1}R_xR_y=R_{y^{-1}xy}$. On the other hand, $L/Z(L)$ has
exponent $2$, therefore $[x,y] \in Z(L)$ for all $x,y \in L$.
\[R_x^{-1} R_y^{-1} R_xR_y = R_x^{-1} R_{y^{-1}xy} = R_x^{-1} R_{x[x,y]} =
R_{[x,y]},\]
hence (i) holds. Using the Moufang identities, one obtains
\begin{eqnarray*}
a[R_y,L_z] &=& z((z^{-1}\cdot ay^{-1})y) \\
&=& z((z^{-1}(ay^{-1}\cdot z)z^{-1})y) \\
&=& z(z^{-1}((ay^{-1}\cdot z)(z^{-1}y))) \\
&=& (ay^{-1}\cdot z)(z^{-1}y) \\
&=& aR_{y^{-1},z}.
\end{eqnarray*}
This shows (ii), since $R_{y^{-1},z} =L_{y,z^{-1}}$ is well known from
\citet[Lemma VII.5.4.]{Bruck}. The inner map $R_{x,y}= R_xR_yR_{xy}^{-1}$
is also known to be a pseudo-automorphism with companion $[x,y]$. Since
$L/Z(L)$ is Abelian, $R_{x,y}$ turns out to be an automorphism. Moreover,
$z^{-1} \, (zR_{y,x}) \in Z(L)$ for all $z \in L$.  We have
\begin{eqnarray*}
[R_x,[R_y,L_z]] &=& [R_x,R_{y^{-1},z}] \\
&=& R_x^{-1} R_{xR_{y^{-1},z}} \\
&=& R_{x^{-1} \, (xR_{y^{-1},z})} 
\end{eqnarray*}
and
\[x^{-1} \, (xR_{y^{-1},z}) = (x,z^{-1},y)^{-1}\]
by \citet[Lemma VII.5.4.]{Bruck}. Using the fact that $L/Z(L)$ elementary
Abelian, we get (iii). \qed

\section{Structure constants of the triality group}

Let $(V,\sigma,\kappa,\alpha)$ be a symplectic cubic space and let us
construct the group with triality $(G,S)$ as in Section
\ref{sec:grconstr}. Write $(L=\sigma^G, \circ)$ for the Moufang loop
associated to $(G,S)$. Since the elements of $H_3=C_G(\sigma)$ are
\[\{ (g_1f_1)^{x_1} \cdots (g_nf_n)^{x_n} \, h_1^{y_1} \cdots h_n^{y_n} \,
(uv)^z \}.\]
The largest normal subgroup $N$ of $G$ contained in $C_G(\sigma)$ contains
the element $uv$. Let $a$ be an arbitrary element of $G$ and let us denote
by $\bar a$ the right action of $a$ on $L=\sigma^G$. (This action is
naturally equivalent with the right action on the right cosets of
$C_G(\sigma)$.) Put $\bar G = G/N = \{\bar a \mid a\in G\}$. Clearly, $\bar
u=\bar v \in Z(\bar G)$. Moreover, the group $A=\langle \bar g_1,\ldots,
\bar g_n, \bar u\rangle$ acts sharply transitively on $L$.

As we saw in Proposition \ref{pr:mltgr}, for the element $x =\sigma^g \in
L$, one has $R_x=\bar \gamma$ with $\gamma= [g,\sigma]^\rho$ and $L_x =
\overline{\gamma^\rho}$. Put $s=\sigma^u$, then $R_s=L_s=\bar u$ and $s \in
Z(L)$. For the element $x=\sigma^{g_1^{x_1} \cdots g_n^{x_n}}$, $R_x=\bar
\gamma$ with
\begin{eqnarray*}
\gamma&=&[g_1^{x_1} \cdots g_n^{x_n},\sigma]^\rho \\
&=& f_n^{-x_n} \cdots f_1^{-x_1} \, (g_1f_1)^{-x_1} \cdots (g_nf_n)^{-x_n}
\\
&=& g_1^{x_1} \cdots g_n^{x_n} \, (\prod_{k; i<j} h_k^{\alpha_{ijk}x_ix_j})
\, u^av^b
\end{eqnarray*}
for some $a,b \in \mathbb Z_2$, that is, 
\[R_x= \bar \gamma = {\bar g}_1^{x_1} \cdots {\bar g}_n^{x_n} \, (\prod_{k;
  i<j} {\bar h}_k^{\alpha_{ijk}x_ix_j}) \, {\bar u}^z\]
for some $z\in \mathbb Z_2$. In particular, the set of right multiplications of
$L$ is
\[ \{ {\bar g}_1^{x_1} \cdots {\bar g}_n^{x_n} \, (\prod_{k; i<j}
   {\bar h}_k^{\alpha_{ijk}x_ix_j}) \, {\bar u}^z \mid x_1,\ldots,x_n,z \in
   \mathbb Z_2 \}.\] 
A consequence of this is $R_x^2 \in \langle \bar u \rangle$, or
equivalently $x^2\in \langle s \rangle$ for all $x\in L$. This means that
$L/\langle s\rangle$ is an elementary Abelian $2$-group and $L$ is a small
Frattini Moufang loop.

We are now able to prove our main result.

\begin{thm}
Let $(V,\sigma,\kappa,\alpha)$ be a symplectic cubic space and let us
construct the group with triality $(G,S)$ using the relations
\eqref{eq:G1}--\eqref{eq:G5} and \eqref{eq:sigma}--\eqref{eq:rho}. Let $L$
be the Moufang loop associated to $(G,S)$. Then, $L$ is a small Frattini
Moufang loops and the symplectic cubic space corresponding to $L$ is $V$.
\end{thm}
\pf It only remained to show that the structure constants of $L$ and $V$
are the same. Let us put $x_i =\sigma^{g_i} \in L$ and define $A=\langle
g_1, \ldots, g_n, u \rangle$ as before. The set $\{g_1,\ldots,g_n\}$ is
independent in $A$, that is, no proper subset of it generates $A$. This
implies that the set $\{x_1, \ldots, x_n\}$ is independent in $L$. With
other words, $\{x_1\langle s \rangle, \ldots, x_n\langle s \rangle\}$ is a
basis for the vector space $L/\langle s \rangle$.

Then, by Proposition \ref{pr:mltgr}(ii), 
\[R_{x_i} = {\bar g}_i \, {(\bar g_i \bar f_i)}^{-2}= {\bar g}_i
\hskip 5mm \mbox{and} \hskip 5mm
L_{x_i} = {\bar f}_i {\bar u}^{\sigma_i}.\]
By Lemma \ref{lm:moufids}, 
\begin{eqnarray*}
R_{x_i}^2 &=& {\bar g}_i^2 = {\bar u}^{\sigma_i}, \\
R_{[x_i,x_j]} &=& [{\bar g}_i,{\bar g}_j] = \bar u^{\kappa_{ij}}, \\
R_{(x_i,x_j,x_k)} &=&  [[{\bar g}_i,\bar f_j],{\bar g}_k] = \bar
u^{\alpha_{ijk}}. 
\end{eqnarray*}
This means
\[ {x_i}^2 = s^{\sigma_i}, \;
{[x_i,x_j]} = s^{\kappa_{ij}}, \; \mbox{ and } \;
{(x_i,x_j,x_k)} = s^{\alpha_{ijk}}, \]
hence the symplectic cubic space of $L$ is indeed $V$. \qed

\bibliographystyle{alpha}

\end{document}